\documentclass[]{amsart}
\usepackage{amsmath, amsthm, amsfonts, amssymb, color,graphicx,mathrsfs}
\newtheorem{thm}{\bf Theorem}

\newtheorem{lem}{\bf Lemma}
\newtheorem{prp}{\bf Proposition}

\theoremstyle{definition}

\newcommand{\bfr}[1]{\mathscr #1}
\definecolor{wco}{rgb}{0.5,0.2,0.3}

 \def\Re{{\rm Re}\;} \def\Im{{\rm Im}\;}
\def\D{\mathbb{D}}
\def\J{\bfr J}
 \def\A{\bfr{A}} \def\B{\bfr{B}} \def\W{\mathbf{W}}
\def\OM{\mathbf{\Omega}}\def\C{\mathbb{C}}
\def\ds{\displaystyle}

\begin{document}

\begin{center}{\bf\large On The Dynamics Of The Rational Family $\displaystyle
\mathbf{f_t(z)=-\frac{t}{4}\frac{(z^{2}-2)^{2}}{z^{2}-1}}$}\bigskip\bigskip

\small\bf Hye Gyong Jang and Norbert Steinmetz
\end{center}

\bigskip\begin{quote} {\small {\bf Abstract.} In this paper we discuss the
dynamics as well as the structure of the parameter space of the
one-parameter family of rational maps $\ds
f_t(z)=-\frac{t}{4}\frac{(z^{2}-2)^{2}}{z^{2}-1}$ with free
critical orbit
$\pm\sqrt{2}\xrightarrow{(2)}0\xrightarrow{(4)}t\xrightarrow{(1)}\cdots$.
In particular it is shown that for any escape parameter $t$ the
boundary of the basin at infinity $\A_t$ is either a Cantor set, a
curve with infinitely many complementary components, or else a
Jordan curve. In the latter case the Julia set is a Sierpi\'nski
curve.


}\end{quote}

\section{\bf Introduction}

The dynamics and the structure of the parameter space of
one-parameter families of rational maps with a single free
critical orbit were studied in
\cite{Bed,BLa,Dev1,Dev2,Dev3,Ste1,Ste2,Ste3} and many other
papers. For any escape parameter $t$ (satisfying
$f_t^{k}(c_t)\to\infty$ for one, hence all free critical values)
in the so-called McMullen family $\ds f_t(z)=z^{m}+t/z^{l}$, the
Julia set of $f_t$ is either a Cantor set, a Sierpi\'nski curve or
else a Cantor set of Jordan curves. In the so-called
Morosawa-Pilgrim family $\ds
f_t(z)=t\Big(1+\frac{\frac{4}{27}z^{3}}{1-z}\Big)$, for any escape
parameter $t$ the Julia set is either a Cantor set, a Sierpi\'nski
curve or else a curve which has either one or else infinitely many
cut-points, hence is not a Sierpi\'nski curve. Similar results
were obtained for the modular family $\ds
f_t(z)=\frac{4t}{27}\frac{(z^2-z+1)^3}{z^2(z-1)^2}$. In this paper
we will consider the degree-four rational family
 $$f_t(z)=-\frac{t}{4}\frac{(z^{2}-2)^{2}}{z^{2}-1}$$
with super-attracting fixed point $\infty$ and free critical orbit
 $$\pm\sqrt{2}\xrightarrow{(2)}0\xrightarrow{(4)}t\xrightarrow{(1)}\cdots.$$
The family $f_t$ is even with respect to $z$, hence the Julia set
$\J_t$ is also even ($\J_t=-\J_t$). Any Fatou component $U$ of $f_t$ is either
symmetric ($-U=U$) and contains $z=0$ or $z=\infty$, or else its
counterpart $-U\ne U$ is also a Fatou component. The family is odd
with respect to $t$, hence $\J_{-t}=-\J_{t}=\J_t$ holds, and the
structure of the parameter plane $\C^*=\C\setminus\{0\}$ also has
this rotational symmetry. The semi-conjugate family
 $\ds\hat f_t(z)=f_t(\sqrt z)^2=\frac{t^2}{16}\frac{(z-2)^4}{(z-1)^2}$
has no longer symmetries in the dynamical plane, but has the same
parameter space. Replacing $t^2$ by $t$ rules the symmetry out.

\medskip{\bf B\"ottcher's function.} For $t\ne 0$ the map $f_t$ has a
super-attracting fixed point at $z=\infty$, the corresponding
super-attracting basin will be denoted by $\A_t$. {B\"ottcher's
function} $\Phi_t$ is the solution to B\"ottcher's functional
equation
 $$\Phi_t\circ f_t(z)=-\frac t4 \Phi_t(z)^2;$$
it is uniquely determined by the normalisation $\Phi_t(z)\sim z$
at $z=\infty$. If $\A_t$ is simply connected, then $\Phi_t$ maps
$\A_t$ conformally onto the disc $|w|>4/|t|.$\medskip

{\bf Bifurcation locus.} The set of $t\in\mathbb{C}^*$ such that
the Julia set of the family $(f_t)$ does not move continuously
over any neighbourhood of $t$ is called {bifurcation locus}, see
McMullen \cite{McM}. By the central theorem in \cite{McMullen00},
small copies of the bifurcation locus of the family $(z,t)\mapsto
z^8+t$ are dense in the {bifurcation locus} (the number $8$ is
explained by the fact that $f_t^2$ has degree $8$ at $z=\pm\sqrt
2$).

\bigskip\includegraphics[viewport=0 0 368 245]{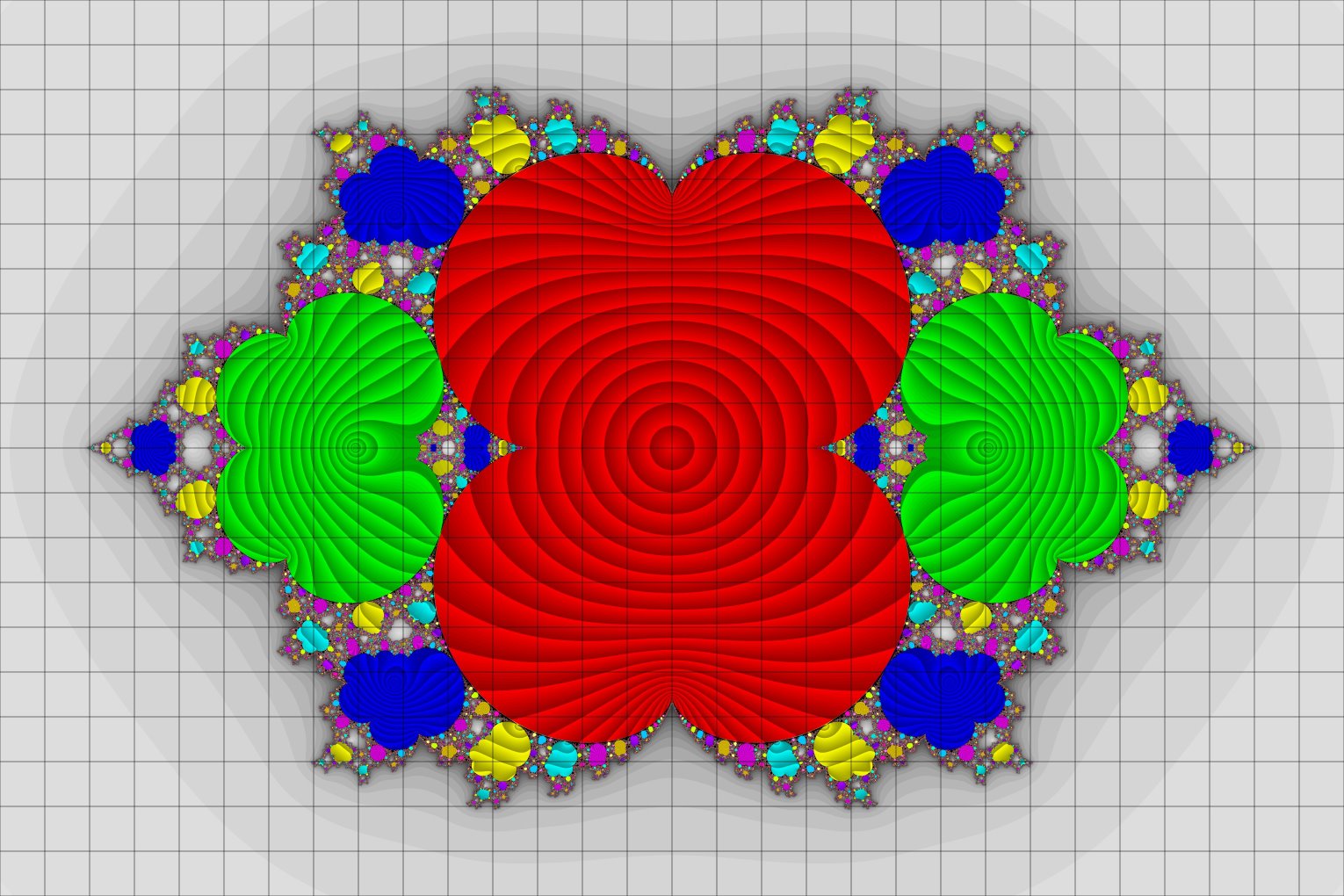}

\medskip

\begin{center}\small{\bf Figure 1.} The parameter space for the family $(f_t)$
($|\Re t|<3$, $|\Im t|<2$). The bifurcation locus is the boundary
of the escape locus (grey). The hyperbolic components are
displayed till period $14$ (red=1, green=2, blue=3, yellow=4,
...).
\end{center}

\medskip{\bf Escape locus.} The set of all $t$ such that
$f^n_t(t)\in \A_t$ for some $n\ge 0$ is called {escape locus}.
More precisely we will write
 $$\begin{array}{rcl}
 \OM_0&=&\{t:t\in \A_t\},\cr
 \OM_n&=&\{t:f_t^n(t)\in
 \A_t,~f_t^{n-1}(t)\notin \A_t\}\quad(n\ge 1), {\rm~and}\cr
\OM_\infty&=&\C^*\setminus \bigcup_{n\ge
0}\OM_n=\{t:f_t^n(t)\notin\A_t {\rm~for~all~} n\};\end{array}$$
the latter set consists of those functions $f_t$ with bounded free
critical orbit (actually it is bounded by $2\sqrt{2}$), while if
$t\in\OM_n$ for some $n\ge 0$, all critical points are attracted
to $\infty$, so $f_t$ is hyperbolic. By the $\lambda$-Lemma due to
Ma\~{n}\'e, Sad, and Sullivan \cite{MSS}, every set $\OM_n$ is
open. The escape set may equally be described by the recursively
defined sequence $(Q_n)$ of rational functions
 $$Q_0(t)=t,~Q_n(t)=f_t(Q_{n-1}(t))=f_t^{n}(t):$$
$t\in \OM_n$ if and only if $Q_n(t)\in \A_t$ and $Q_{n-1}(t)\notin
\A_t$ for $n>0$, and $Q_0(t)=t\in\A_t$ for $n=0$.\medskip

{\bf Hyperbolic components.} Let $\W_n$ be the set of parameters
$t$ such that $f_t$ has a finite (super)-attracting cycle of exact
period $n$. The components of $\W_n$ are called {\it hyperbolic}.
By the Implicit Function Theorem the hyperbolic components are
open.

\section{\bf Cantor And Sierpi\'nski Curve Julia Sets}

The map $f_t$ has a super-attracting fixed point at $\infty$,
hence a super-attracting immediate basin $\A_t$ about infinity;
$\A_t$ is symmetric with respect to the origin, and
$f_t:\A_t\longrightarrow\A_t$ has degree two or four. By $T$ we
denote the component of $f^{-1}_t(\A_t)$ that contains the pole
$z=1$, so $-T$ contains the pole $-1$. Since
$0\xrightarrow{(4)}t$, the point $t$ belongs to $\A_t$ if and only
if $\A_t$ is completely invariant and contains all critical
points, hence coincides with the Fatou set. Since $f_t$ is
hyperbolic, we have (see any of \cite{Bea,Mil,Ste})

\begin{thm}\label{T3.1} The Julia set $\J_t$ is a Cantor set if and
only if  $t\in \A_t$, i.e.\ if and only if $t\in\OM_0$.
\end{thm}

 The set $\OM_0$ is also called {\it Cantor locus.}

\begin{thm}\label{T3.2}For~ $t\notin \OM_0$, the immediate basin
$\A_t$ is simply connected, and the domains $-T$ and $T$ are
distinct.\end{thm}

{\bf Proof.} Since $t\notin \A_t$ implies $0\notin \A_t$, the
fixed point at $\infty$ is the only critical point and value in
$\A_t$, hence (again by any of \cite{Bea,Mil,Ste}) $\A_t$ is
simply connected. Assuming $-T=T$ and noting that $T$ contains no
critical point, it follows that $f_t:T\longrightarrow \A_t$ is an
{\it unramified} covering of degree two. This, however, is
impossible since $\A_t$ is simply connected.~$\square$

\begin{prp}\label{P3.4}Suppose
$\overline{T}\cap\overline{\A_t}=\emptyset$, and let $D$ be the
component of $\overline{\C}\setminus\overline{\A_t}$ containing
$T$. Then $f_t(D)=\overline{\C}$, $D$ is symmetric with respect to
the origin and contains $z=0$, $-T$ and $f_t^{-1}(D)$ {\rm (of
course, $T$ and $-T$ may be interchanged)}.
\end{prp}

{\bf Proof.} From $\partial D\subset\partial \A_t$ and $\partial
T\subset D$ follows
 $$\partial f_t(D)\subset f_t(\partial D)\subset f_t(\partial \A_t)=\partial
 \A_t=f_t(\partial T)\subset f_t(D),$$
thus $f_t(D)=\overline{\C}$. Also from $t\in f_t(D)$ follows $0\in
D$, hence $D$ is symmetric with respect to the origin and contains
also $-T=-T$ and the points $\pm\sqrt{2}$ (the pre-images of $0\in
f_t(D)$). Since any component $D'$ of $f_t^{-1}(D)$ is contained
in $\overline{\C}\setminus\overline{\A_t}$ and contains one of the
points $\pm\sqrt{2}$, it is contained in $D$, which thus is
backward invariant.~$\square$

\begin{thm}\label{P3.5}Suppose $t\in \OM_n$ for some $n>0$. Then the following are equivalent:

\begin{itemize}\item[(a)] $\overline{T}\cap\overline{\A_t}=\emptyset$
$[$equivalently $\overline{-T}\cap\overline{\A_t}=\emptyset]$.
\item[(b)] $\partial\A_t$ is a Jordan curve. \item[(c)] The Julia
set $\J_t$ is a Sierpi\'nski curve.\end{itemize}\end{thm}

{\bf Proof.} The Julia set is connected and locally connected by
the hypothesis ``$t\in \OM_n$ for some $n>0$''. Thus
``(a)$\Rightarrow$(b)'' follows from Proposition \ref{P3.4} and
Morosawa's Lemma \cite{Mor}(\footnote{{\bf Morosawa's Lemma.} {\it
Any simply connected Fatou fixed domain $\A$ of some
sub-hyperbolic rational map $R$ is bounded by a Jordan curve,
provided there exists some domain $D$ complementary to
$\overline{\A}$ that contains a Fatou component $T$ and also its
pre-image $R^{-1}(T)$.}}). The proof of ``(b)$\Rightarrow$(c)'' is
a consequence of Lemma 5 in \cite{Ste1}(\footnote{{\bf Lemma 5.}
{\it Let $R$ be rational and $D$ be a Jordan domain such that
$\partial D$ contains no critical value of $R$. Then any two
different components of $R^{-1}(D)$ $($are also Jordan domains
and$)$ have disjoint closures.}}), applied to $D=\A_t$, any two
pre-images $D_1$ and $D_2$ of $D$ of any order (including the case
$D_1=D$), and $R=f_t^n$ for $n$ chosen in such a way that
$f^n(D_\nu)=D$ $(\nu=1,2)$. Finally, ``(a)'' is part of the notion
``Sierpi\'nski curve'', namely $\overline{\A_t}\cap
\overline{(-T)}=\overline{\A_t}\cap\overline T=\emptyset,$ hence
``(c)$\Rightarrow$(a)''.~$\square$

\bigskip\includegraphics[viewport=0 0 368 123]{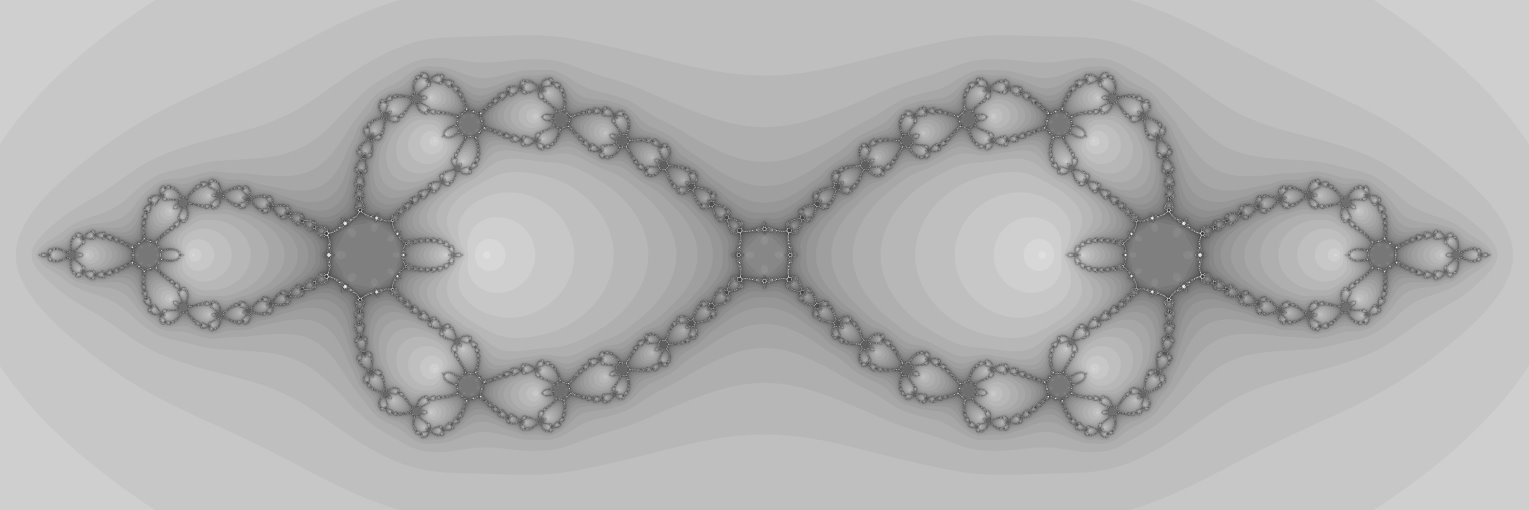}

\bigskip\includegraphics[viewport=0 0 368 123]{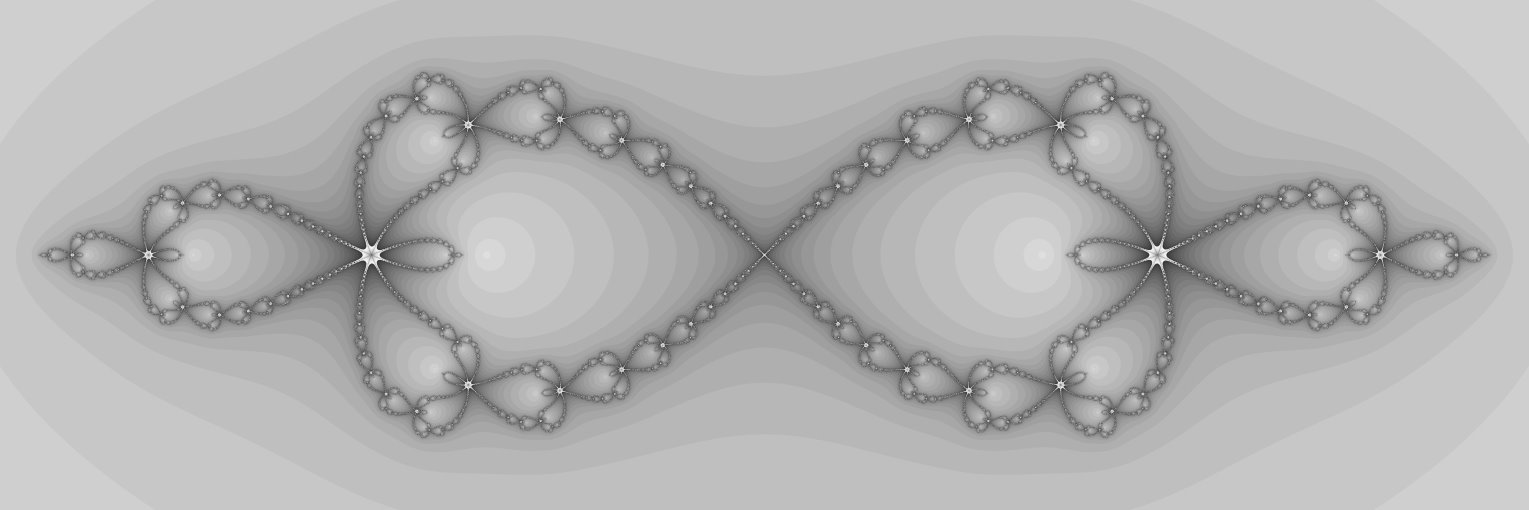}

\medskip
\begin{center}\small {\bf Figure 2.} The Julia set for $t=\sqrt{4+2\sqrt{2}}-0.0001$
is a Sierpi\'nski curve, while for $t=\sqrt{4+2\sqrt{2}}$ (a
Misiurewicz point satisfying $f_t(t)=f_t^2(t)$) it is
not.\end{center}

\begin{thm}If $t$ is a {\it Misiurewicz point,} i.e. if $t\in\J_t$
is pre-periodic, then $\overline{\C}\setminus\overline{\A_t}$
consists of countably many components.\end{thm}

{\bf Proof.} $f_t$ is sub-hyperbolic, hence the Julia set is
connected and locally connected, but $\J_t$ has
self-intersections: one at $z=0$, and hence infinitely many at the
pre-images of $z=0$ under the iterates $f^n_t$; in particular,
$\overline{\C}\setminus\overline{\A_t}$ consists of countably many
components.~$\square$

\begin{thm}For $t\in\bigcup_{n\ge 1}\OM_n$ the basin $\A_t$ is either a Jordan domain
$($and the Julia set $\J_t$ is a Sierpi\'nski curve$)$ or else
$\partial\A_t$ has infinitely many complementary
components.\end{thm}

{\bf Proof.} The complementary components of
$\overline\C\setminus\overline{\A_t}$ are Jordan
domains(\footnote{This is true in general: If a domain $\A$ is
bounded by a curve, then every complementary component of
$\overline\A$ is a Jordan domain.}). Let $D$ be the component
containing $T$. Then $-D$ contains $-T$, and $\A_t$ is a Jordan
domain if and only if $D=-D$, hence $0\in D$. We thus may assume
$-D\cap D=\emptyset.$ Since any component of
$\overline\C\setminus\overline{\A_t}$ different from $\pm D$ is
mapped eventually onto one of $\pm D$,
$\overline\C\setminus\overline{\A_t}$ consists of countably many
components.~$\square$

\medskip{\bf Remark.} $f_t(D)$ contains $\A_t$
and also part of $\partial\A_t$, namely $f_t(D\cap\partial T)$,
thus $f_t(D)$ covers the component of $\overline\C\setminus
f_t(\overline{\A_t}\cap\overline T)$ that contains $\infty$. If
$\A_t$ is not a Jordan domain, then $f_t(D)$ cannot contain the
point $t$, since otherwise $D$ would contain the origin. In
particular, $\overline\C\setminus f_t(\overline{\A_t}\cap\overline
T)$ cannot be connected, so that $\overline{\A_t}\cap\overline T$
must contain a continuum, hence also a Jordan arc. All graphical
experiments, however, support the conjecture that for
$t\in\bigcup_{n\ge 1}\OM_n,$ the Julia set $\J_t$ is a
Sierpi\'nski curve.

\section{\bf The Escape Locus}
Recall that in a neighbourhood of $\infty$, B\"{o}ttcher's
function is given by
\begin{equation}\label{Boettlimit}
\Phi_t(z)=\lim_{k\to\infty}\sqrt[2^{k}]{-f_t^{k}(z)(4/t)^{1+2+\cdots
+2^{k-1}}}=\frac 4t\lim_{k\to\infty}\sqrt[2^{k}]{-\frac t4f_t^{k}(z)}.
\end{equation}
It admits unrestricted analytic continuation throughout $\A_t$
with all bran\-ches satisfying
$\lim\limits_{z\to\partial\A_t}|\Phi_t(z)|=4/|t|$ and
$|\Phi_t(z)|>4/|t|$. If $\A_t$ is simply connected (i.e.\ if
$t\notin \A_t$, equivalently $t\notin\OM_0$), then
(\ref{Boettlimit}) holds throughout $\A_t$.

\medskip{\bf Remark.} The boundary of the Mandelbrot set $\mathcal M_2$ is the
bifurcation locus of the quadratic family $f_t(z)=z^2+t$. Douady
and Hubbard \cite{DH82} proved that the Mandelbrot set is
connected by showing that $\mathbf{\Phi}(t)=\Phi_t(t)$ maps the
complement of $\mathcal M_2$ conformally onto $\triangle_1$. Here
$\Phi_t$ is B\"ottcher's function for the quadratic map $f_t$;
like in our case it depends analytically on $t$. This method has
been used by various authors to prove similar results. Since for
$t\notin\mathcal M_2$ the corresponding attracting basin $\A_t$ is
infinitely connected, the main problem consists in continuing
B\"ottcher's function to some simply connected subdomain of $\A_t$
containing the critical value $t$ and the fixed point $\infty$.

The same difficulty occurs here. It is overcome by constructing an
exhaustion $(D_n)$ of $\A_t$, such that $f_t:D_n\longrightarrow
D_{n-1}$ has degree two as long as $0\not\in D_n$; then $t\in D_m$
and $0\in D_{m+1}$ for the first time, and $D=D_m$ is a simply
connected sub-domain of $\A_t$ we are looking for. So we may
define
 $$E_0(t)=-\frac t4\Phi_t(t)\quad(t\in\OM_0)\quad{\rm with~} E_0(t)\sim -\frac{t^2}4\quad(t\to\infty).$$
Since for $t\not\in\OM_0$ the basin $\A_t$ is simply connected,
there is no difficulty to define
$$E_n(t)=-\frac t4\Phi_t(Q_n(t))
 \quad(t\in\OM_n,~n\ge 1).$$
The function $E_n$ $(n\ge 0)$ is analytic on every component
$\Omega$ of $\OM_n$ and satisfies $|E_n(t)|\to 1$ as
$t\to\partial\Omega$, and $|E_n(t)|>1$ on $\Omega$, hence $\Omega$
contains a pole of $Q_n$. From (\ref{Boettlimit}) follows
 $$E_n(t)=\lim_{k\to\infty}\sqrt[2^k]{-\frac t4Q_{n+k}(t)}\quad(n\ge 0).$$

\begin{thm}\label{T4.2}The components $\Omega$ of $\OM_n$ with $n\ge 1$ are simply
connected, and the restriction $E_n|_\Omega$ is a conformal map
$\Omega\longrightarrow\triangle_1$. The finite poles of $Q_n$ are
simple and distributed over the components of $\OM_n$ in such a
way that each component contains exactly one pole of
$Q_n$.\end{thm}

{\bf Remark.} For $n\ge 1$ fixed, each component of $\OM_n$
contains exactly one solution to the equation $Q_{n-1}(t)^2=1.$ In
particular, the number of these components is $2\deg
Q_{n-1}(t)=2(4^{n}-1)/3$.

\medskip{\bf Proof.} We know that $E_n|_\Omega$ is a proper map
$\Omega\longrightarrow\triangle_1$. The method of Roesch
\cite{Roesch} enables us to show that it is a local
homeo\-morphism, hence is a conformal map since $\triangle_1$ is
simply connected. To this end take any $t_0\in \Omega$ and choose
$\varepsilon>0$ such that $|t-t_0|<3\varepsilon$ belongs to the
Fatou component of $f_{t_0}$ containing $t_0$. Furthermore let
$\eta_t:\overline{\C}\longrightarrow \overline{\C}$ be any
diffeomorphism such that $\eta_t(w)$ depends analytically on $t$
for $|t-t_0|<\varepsilon$ and satisfies
 $\eta_t(w)=w$ on $|w-t_0|\ge 3\varepsilon$ and $\eta_t(w)=w+(t-t_0)$ on
 $|w-t_0|<\varepsilon.$ Then $g_t=\eta_t\circ f_{t_0}:\overline{\C}\longrightarrow
\overline{\C}$ is a quasiregular map that coincides with $f_{t_0}$
on $\overline\C\setminus f^{-1}_{t_0}(\{w:|w-t_0|\le
3\varepsilon\})$, and is analytic on $\overline{\C}\setminus
f_{t_0}^{-1}(A)$, with
$A=\{w:\varepsilon\leq|w-t_0|\leq3\varepsilon\}$. Since the sets
$f_{t_0}^{-n}(A)$ ($n=1,2,\ldots$) are mutually disjoint, $g_t$ is
quasiconformally conjugate to some rational function
 $$R_t=h_t\circ g_t\circ h_t^{-1}$$
by Shishikura's qc-Lemma \cite{Shi}. If the quasiconformal map
$h_t:\overline\C\longrightarrow\overline\C$ is normalised to fix
the critical points $\sqrt{2}$, $0$ and $\infty$, it depends
analytically on the parameter $t$. We set $h_t(-\sqrt{2})=z_0 $,
$h_t(1)=z_1$, and $h_t(-1)=z_{-1}$ to obtain
 $$R_t(z)=a(t)\frac{(z-\sqrt{2})^{2}(z-z_0)^{2}}{(z-z_1)(z-z_{-1})}.$$
Since $h_t(0)=0$, $R_t$ has a critical point of order three at
$z=0$, and solving $R_t'(0)=R_t''(0)=R_t'''(0)=0$ yields $z_1=1$,
$z_{-1}=-1$, and $z_0=-\sqrt{2}$, thus
 $$R_t(z)=a(t)\frac{(z^{2}-2)^{2}}{z^{2}-1}\quad{\rm and}\quad R_t(0)=-4a(t).$$
From $R_t\circ h_t=h_t\circ\eta _t\circ f_{t_0}$ and $h_t(0)=0$,
however, follows
 $$R_t(0)=h_t\circ\eta_t\circ f_{t_0}(0)=h_t\circ\eta_t(t_0)=h_t(t),$$
hence $R_t=f_{h_t(t)}.$ In $\A_{h_t(t)}$ we have
$$\begin{array}{rcl}\frac{t_0}{h_t(t)}\Phi_{t_0}\circ
h_t^{-1}\circ f_{h_t(t)} &=&\frac{t_0}{h_t(t)}\Phi_{t_0}\circ
g_t\circ h_t^{-1}\cr
&=&\frac{t_0}{h_t(t)}\Phi_{t_0}\circ\eta_t\circ f_{t_0}\circ
h_t^{-1}\cr &=& \frac{t_0}{h_t(t)}\Phi_{t_0}\circ f_{t_0}\circ
h_t^{-1}\cr &=&\frac{t_0}{h_t(t)}(-\frac{t_0}{4})(\Phi_{t_0}\circ
h_t^{-1})^{2}\cr&=& -\frac{h_t(t)}{4}(\frac{t_0}{h_t(t)}
\Phi_{t_0}\circ
 h_t^{-1})^{2},\end{array}$$
hence
 $$\Phi_{h_t(t)}=\frac{t_0}{h_t(t)}\Phi_{t_0}\circ h_t^{-1} {\rm ~and~}
 E_n(h_t(t))=-\frac{t_0}{4}\Phi_{t_0}(f_{t_0}^{n}(t)).$$
The right hand side is locally univalent on
$\A_{t_0}\setminus\{\infty\}$, and so is $t\mapsto E_n(h_t(t))$
(note that $t\mapsto h_t(t)$ is analytic on some neighbourhood of
$t_0$). Thus $E_n|_\Omega$ is a conformal map
$\Omega\longrightarrow\triangle_1$. In particular, $Q_n$ has only
simple poles, exactly one in $\Omega$.~$\square$

\begin{prp}\label{P4.1}$\OM_0$ contains the punctured disc
$\{t:3\le|t|<\infty\}$. For~ $|t|\ge 3$ the disc
$\triangle_3=\{z:|z|>3\}$ $($including $\infty)$ is invariant
under $f_t$ and belongs to $\A_t$.\end{prp}

{\bf Proof.} For $|z|\ge 3$ and $|t|\ge 3$ easily follows
 $$|f_t(z)|\ge\frac 34\frac{(|z|^2-2)^2}{|z|^2+1}\ge \frac 34\cdot\frac{49}{10}>3.$$
Thus we have $t\in\OM_0$, the Julia set is totally disconnected,
and $\A_t$ contains $\triangle_3$ (one can even replace the number
$3$ by the smaller number $2\sqrt{2}$).~$\square$

\begin{thm}The {\em Cantor locus} $\OM_0$ is connected, and $\Omega_0=\OM_0\cup\{\infty\}$ is
mapped conformally onto the disc $\triangle_1=\{w:|w|>1\}$ by any
branch of $\sqrt{E_0}$.\end{thm}

{\bf Proof.} Let $\Omega$ be any (in $\C$) bounded component of
$\OM_0$. Then $E_0$ is holomorphic on $\Omega$ and satisfies
$|E_0(t)|\to 1$ as $t\to\partial\Omega$. The maximum principle
yields $|E_0(t)|<1$ in contrast to $|E_0(t)|>1$. Thus $\OM_0$ is
connected and contains $3\le|t|<\infty$. Now repetition of the
proof of Theorem \ref{T4.2} shows that $E_0$ is locally univalent
on $\OM_0$, thus has no critical points in $\Omega_0$ except at
$t=\infty$. Since $E_0:\Omega_0\longrightarrow\triangle_1$ has
degree two, the Riemann-Hurwitz formula shows that $\Omega_0$ is
simply connected and $\sqrt{E_0(t)}=\frac i2t+\cdots$ is a
conformal map of $\Omega_0$ onto $\triangle_1$.~$\square$

\section{\bf Kernel Convergence}
 The concept of kernel convergence in the sense of
Carath\'{e}odory may be described as follows (see Pommerenke
\cite{Pommi}): Let $(D_n)$ be any sequence of domains, each
containing some base point $z_0$. The {\it kernel} ker$(D_n)=K$ of
$(D_n)$ with respect to $z_0$ then is the union of all simply
connected domains $D$, such that $z_0\in D$ and
$\overline{D}\subset D_n$ for $n\geq n_0(D)$. If no such $D$
exists, we set $K=\{z_0\}$. The sequence $(D_n)$ is said to {\it
converge} to $K$ in the sense of Carath\'{e}odory, if $K$ is also
the kernel of every sub-sequence $(D_{n_k})$. If all domains $D_n$
are simply connected and if $f_n$ denotes the conformal map of the
unit disc $\D$ onto $D_n$, normalised by $f_n(0)=z_0$ and
$f'_n(0)>0$, then $D_n\longrightarrow K\neq\{z_0\}$ in the sense
of Carath\'{e}odory is equivalent to $f_n\to f$ locally uniformly,
where $f$ is the normalised conformal map $\D\longrightarrow K;$
for $K=\{z_0\}$ the sequence $(f_n)$ tends to the constant $z_0$.

The dynamical spheres for the functions $Q_n(t)=f^n_t(t)$  and the
parameter space for the family $(f_t)$ are intimately connected
via kernel convergence. This was first observed by N.\ Busse
\cite{Buss} in case of the quadratic family $f_t(z)=z^2+t$. We
start with simple properties of $Q_n$.

\begin{lem}$Q_n$ satisfies
$Q_n(t)\sim a_n t^{2^{n}-1}$ as $t\to\infty$ with leading
coefficient $a_n=-1/4^{2^{n}-1}$, and has degree
$(4^{n+1}-1)/3.$\end{lem}

{\bf Proof.} From the recursion formula $Q_n(t)=f_t(Q_{n-1}(t))$
follows
 $$a_nt^{\alpha_n}=-\frac14a_{n-1}^2t^{1+2\alpha_{n-1}},$$
hence $2^{-n}(\alpha_n+1)=2^{-(n-1)}(\alpha_{n-1}+1)$ and
 $$\alpha_n+1=2^{n}(\alpha_0+1)=2^{n+1}$$
on one hand, and $a_n=-\frac 14 a_{n-1}^2$, hence
$a_n=(-1/4)^{2^{n}-1}$, on the other. Again from the recursion
formula and $Q_n(0)=0$ follows
 $$q_n=\deg Q_n=4q_{n-1}+1,$$ hence
$4^{-n}(q_n+1/3)=4^{-n-1}(q_{n-1}+1/3)$ and
$q_n+1/3=4^n(q_0+1/3)=4^{n+1}/3.$~$\square$\medskip

The rational function $Q_n$ has a super-attracting fixed point at
$t=\infty$ with super-attracting basin $\B_n$ about $\infty$, and
a corresponding B\"{o}ttcher function $\Xi_n$ satisfying
$$\Xi_n(Q_n(t))=a_{n}\Xi_n(t)^{\alpha_{n}}\quad(\Xi_n(t)\thicksim t
 {\rm~as~} t\to\infty),$$
with $\alpha_n=2^{n+1}-1$ and $a_n=-4^{-(2^n-1)}$. Since
$|f_t(z)|>3$ for $|t|\ge 3$ and $|z|\ge 3$, $\B_n$ contains the
closed disc $\overline{\triangle_3}$, hence the kernel of the
sequence $(\B_n)$ with respect to $t=\infty$ contains
$\triangle_3$, and
 $$\Xi_n(t)=\lim_{k\to\infty}\sqrt[\alpha_n^k]{Q_n^k(t)/a_n^{1+\alpha_n+\cdots+\alpha_n^{k-1}}}$$
holds (at least) on $\triangle_3$.

\begin{thm} \label{T5.7}The sequence $(\Xi_n)$ tends to the branch
of $\sqrt{-4E_0}$ satisfying $\sqrt{-4E_0(t)}\sim t$ as
$t\to\infty$, locally uniformly on $\OM_0\cup\{\infty\}$, while
the sequence $(\B_n)$ tends to its kernel $\OM_0\cup\{\infty\}$
with respect to $\infty$.\end{thm}

\medskip{\bf Proof.} Let $D$ be any simply connected domain such that
$\overline{D}\subset \OM_0\cup\{\infty\}$ and $\infty\in D$. Since
$Q_n(t)\to\infty$, uniformly on $D$, we have $|Q_n(t)|>3$ on $D$,
hence $Q_n(D)\subset \B_n$ $(n\ge n_0)$, and $\Xi_n$ is defined on
$Q_n(D)$. Thus also
 $$\sqrt{E_0(t)}=\lim_{n\to\infty}\sqrt[2^{n+1}]{-\frac t4Q_n(t)}=\lim_{n\to\infty}\Psi_n(t)$$
is defined on $D$. From
$$|f_t(z)|\le \frac34\frac{(|z|^2+2)^2}{|z|^2-1}<
 \max\{25,|z|^2\}\quad(|t|=3,~|z|\ge 3)$$
easily follows $|Q_n(t)|\le 5^{2^n}$ $(|t|=3)$ by induction on
$n$, hence
$$|Q_n(t)|\le 5^{2^n}\frac{|t|^{\alpha_n}}{3^{\alpha_n}}\le (5/3)^{2^n}|t|^{\alpha_n}\quad(|t|\ge
 3)$$
by the maximum principle. Also by induction, this time on $k$,
follows
$$|Q_n^k(t)|/|a_n|^{1+\alpha_n+\cdots+\alpha_n^{k-1}}\le
\big(20/3\big)^{2^n(1+\alpha_n+\cdots+\alpha_n^{k-1})}|t|^{\alpha_n^k}$$
hence $|\Xi_n(t)|=\lim\limits_{k\to\infty}\sqrt[\alpha_n^k]
{|Q_n^k(t)|/|a_n|^{1+\alpha_n+\cdots+\alpha_n^{k-1}}}\le 20$ 
on $|t|=3$. The maximum principle then yields $|\Xi_n(t)-t|\le 23$
on $|t|>3$, hence
 $$a_n\Xi_n(t)^{\alpha_{n}}=\Xi_n(Q_n(t))=Q_n(t)+O(1)\quad(|t|>3)$$
with $O(1)$-term independent of $n$, and
 $$\begin{array}{rcl}a_n^{1/\alpha_{n}}\Xi_n(t)&=&\Xi_n^{1/\alpha_n}(Q_n(t))
 =Q_n(t)^{1/\alpha_n}(1+O(\alpha_n^{-1}))\cr
 &=&\Psi_n(t)\big(-4\Psi_n(t)/t\big)^{1/\alpha_n}(1+O(\alpha_n^{-1}))
 \to\sqrt{E_0(t)}\end{array}$$
as $n\to\infty$. By choosing the roots appropriately we obtain
$\Xi_n(t)\to \sqrt{-4E_0(t)}$ and $\OM_0\cup\{\infty\}\subset$
ker$(\B_n).$

Conversely, let $K_1$ be the kernel of any sub-sequence
$(\B_{n_k}),$ and $D$ be any simply connected domain containing
$\infty$ and such that $\overline{D}\subset K_1$. Then
$\overline{D}\subset \B_{n_k}$ $(k>k_0)$, and $\Xi{n_k}$ is
defined on $D$. From
$$|\Xi{n_k}(t)|>4^{\frac{2^{n_k}-1}{2^{{n_k}+1}-2}}=2$$
follows that the sequence $(\Xi{n_k})_{k\geq k_0}$ is normal on
$D$. By Vitali's Theorem on the pointwise convergence of normal
sequences, $(\Xi{n_k})$ converges to $\Xi$, locally uniformly on
$D$, and
 $$\Xi(t)=\sqrt{-4E_0(t)}$$
holds on $\OM_0\cap D$. We claim that
$D\subset\OM_0\cup\{\infty\}$. For otherwise there would be some
$a\in D\cap\partial \OM_0$, hence $|\Xi(a)|=\lim\limits_{t\to
a}\sqrt{|-4E_0(t)|}=2$, this contradicting $|\Xi(t)|\ge 2$, hence
$|\Xi(t)|>2$ on $D$ by the minimum principle. We thus have proved
$D\subset \OM_0\cup\{\infty\}$ and
$K_1\subset\OM_0\cup\{\infty\}$, and so $\OM_0\cup\{\infty\}$ is
the kernel of $(\B_n)$.~$\square$

\bigskip\includegraphics[viewport=0 0 178 80]{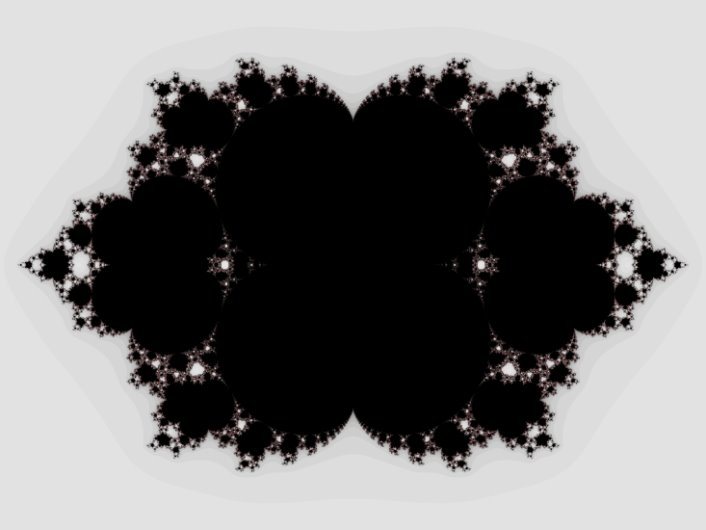}\includegraphics[viewport=0 0 180 133]{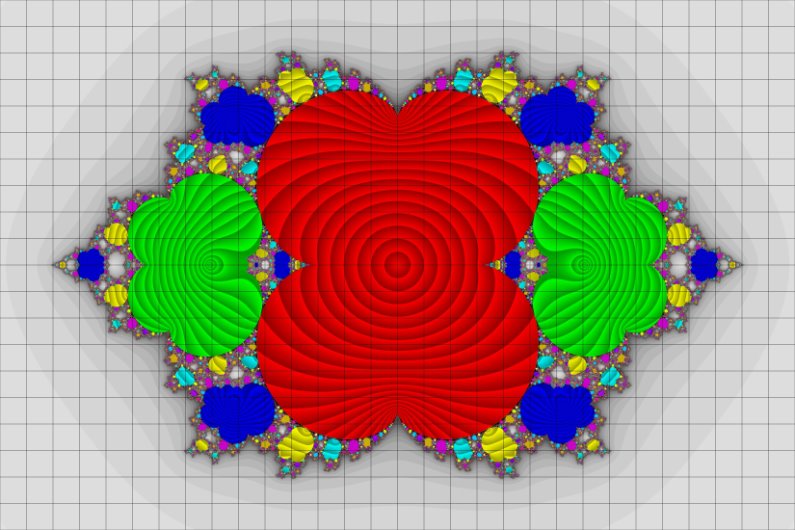}
\medskip

\begin{center}\small{\bf Figure 3.} The dynamical plane for $Q_2$ (left) and the parameter space.\end{center}

\begin{thm} \label{T5.8}Let $t_n$  be any solution to the equation
$Q_{n-1}(t)^2=1$ for some $n\ge 1$, and let $B_k$ $(k\ge n)$
denote the Fatou component of $Q_k$ containing  $t_n$. Then the
sequence $(B_k)$ tends to its kernel $K$ with respected to the
point $t_n$, and $K$ is the connected component of $\OM_n$ that
contains $t_n$.\end{thm}

\section{\bf Hyperbolic Components}
 If $f_t$ has a finite (super)-attracting $n$-cycle
$\{z_0,z_1,\ldots ,z_{n-1}\}$ with associated Fatou cycle $\{U_0
,U_1 ,\ldots ,U_{n-1}\}$, then $\bigcup_{\nu=0}^{n-1}U_\nu$
contains at least one critical point, and from
$\pm\sqrt{2}\xrightarrow{(2)}0\xrightarrow{(4)}t$ follows that the
cycle contains $0$ and at least one of $\pm\sqrt{2}$. We may
assume $0\in U_1$, hence $\sqrt{2}\in U_0$, say.
The domain $U_1$ is symmetric with respect to the origin
$(-U_1=U_1)$. We note that $U_0=U_1$ if $n=1$. The map
$t\mapsto\lambda_t=(f_t^{n})'(z_0)$ is called {\it multiplier
map}, it satisfies $|\lambda_t|\to 1$ as $t\to\partial
W\setminus\{0\}$ on any component $W$ of $\W_n$.

\begin{thm}\label{T6.1}$\W_1\cup\{0\}$ is a simply connected domain;
it is mapped onto the unit disc $\D$ by the multiplier map with
mapping degree four $(\lambda_t\sim t^4$ as $t\to 0)$.\end{thm}

{\bf Proof.} Since none of the maps $f_t$ ($t\in \W_1$) has a
finite super-attracting fixed point (none of the critical points
$\pm\sqrt{2}$ and $0$ is a fixed point), the multiplier map has no
zeros on $\W_1$. As $t\to 0$, $f_t(z)$ has a fixed point $z_t
=t+O(t^{3})$ with multiplier $\lambda_t =t^{4}+O(t^{6}),$ thus
$\W_1$ consists of a single component about the origin $t=0$. To
prove that $\lambda_t$ has no critical points on $\W_1$ we set
 $$\ds f_t(z)=-\frac t4g(z),~g(z)=\frac{(z^{2}-2)^{2}}{z^{2}-1} {\rm ~and~}h(z)=\frac z{g(z)}.$$
Writing $z$ and $\lambda$ instead of $z_t$ and $\lambda_t$ we
obtain from $z=f_t(z)$ and $\lambda=f_t'(z)$
$$-\frac t4=h(z),\quad\lambda=(hg')(z),\quad
\dot{\lambda}=(hg')'(z)\,\dot z,\quad
 -\frac14=h'(z)\,\dot z$$
(the dot means differentiation with respect to $t$), and finally
 $$\dot{\lambda}=-\frac{(hg')'(z)}{4 h'(z)}=
 -\frac{(z^2-2)z^3(3z^2-4)}{4(z^4+3z^2-2)(z^2-1)^2}.$$
The right hand side has the zeros $0,$ $\pm\sqrt{2}$ and
$\pm\sqrt{4/3}$. Now $z=0$ and $\pm\sqrt{2}$ are not fixed points
of any $f_t$, and for $z=\pm\sqrt{4/3}$ we obtain
$\ds\lambda=zg'(z)/g(z)=-16,$ hence $t$ is not in $\W_1$. Thus on
$\W_1\cup\{0\}$ the multiplier map $t\mapsto\lambda_t$ has degree
four with a single critical point of order three at $t=0$, and by
the Riemann-Hurwitz formula $\W_1\cup\{0\}$ is simply
connected.~$\square$

 \hspace*{13mm}\includegraphics[viewport=0 0 368 245]{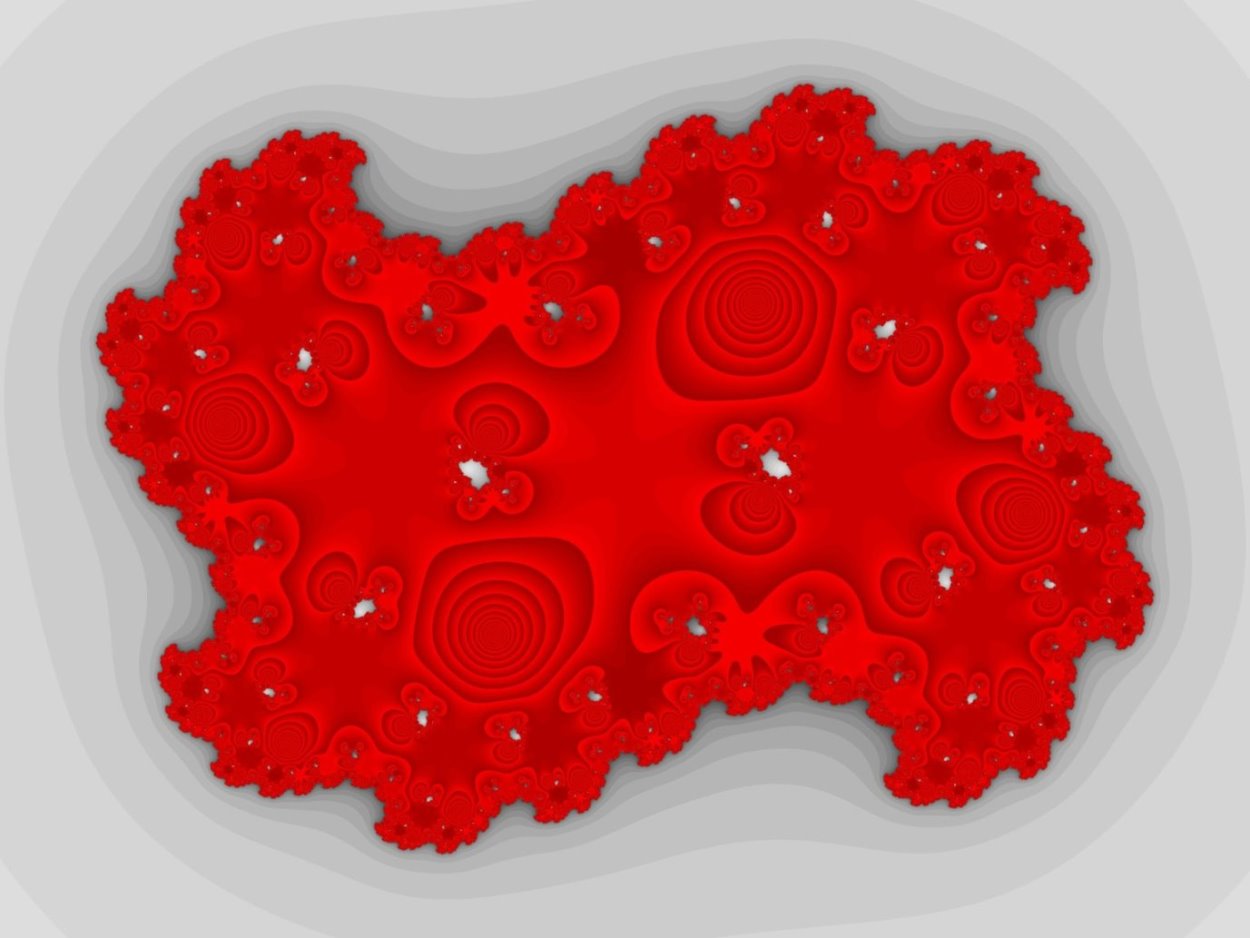}
\medskip

\begin{center}\small {\bf Figure 4.} $f_{0.4+i\,1.3}$ has an attracting fixed
point with basin $S_t$ of infinite connectivity.\end{center}

\begin{thm} \label{T6.3}For\, $t\in \W_1$, the Fatou set of $f_t$
consists of the simply connected domain $\A_t$, its simply
connected pre-images of any order, and the completely invariant
attracting basin $S_t$ about the finite attracting fixed point
$z_t$ of $f_t$. The basin $S_t$ is infinitely connected, while
$\A_t$ is bounded by a Jordan curve.
\end{thm}

{\bf Proof.} The basin $S_t$ contains the critical point $0$ and
the critical value $t=f_t(0)$, hence $f_t:S_t\longrightarrow S_t$
has degree four, and $S_t$ is completely invariant. It is
infinitely connected by the Riemann-Hurwitz formula since it
contains five critical points $-\sqrt{2}, \sqrt{2},0,0,0$. Since
$f_t$ is hyperbolic and $\A_t$ and its pre-images are simply
connected, these components are bounded by curves. To prove that
these curves are Jordan curves it suffices to prove that $\A_t$ is
a Jordan domain.

Our first proof relies on Morosawa's Lemma. Let $D$ denote any
component of $\overline{\C}\setminus\A_t$. Then $D$ contains the
domain $S_t$, hence is uniquely determined. Since
$\overline{\A_t}$ is forward invariant,
$D=\overline{\C}\setminus\overline{\A_t}$ is a backward invariant
domain, and so $\partial \A_t$ is a Jordan curve by Morosawa's
Lemma.

Our second proof will in addition show that the boundary of the
Cantor locus $\OM_0$ is contained in $0.79\le|t|\le 2\sqrt{2}$.
The conjugate
 $$F_t(z)=1/f_t(1/z)=\frac{-4z^2(1-z^2)}{t(1-2z^2)^2}$$
has a super-attracting fixed point at the origin. We have
 $$|F_t(z)|>\frac{4\rho^2(1-\rho^2)}{|t|(1+2\rho^2)^2}\quad(|z|=\rho<1),$$
hence $|F_t(z)|>\rho$ if $|t|\ds\le
\frac{4\rho(1-\rho^2)}{(1+2\rho^2)^2}.$ Then the component $D_t$
of $F_t^{-1}(\D_\rho)$ (with $\D_\rho=\{z:|z|<\rho\}$) that
contains the origin is compactly contained in $\D_\rho,$ and
 $$F_t:D_t\longrightarrow\D_\rho$$
is a polynomial-like mapping (note that $F_t(\pm 1)=0$ and $\deg
F_t|_{D_t}=2$), hence $F_t$ is hybrid equivalent to $z\mapsto
z^2$. The Julia set of $F_t$ thus is a quasicircle, and so is the
boundary curve of $\A_t$ for $0<|t|<\tau$; the value $\tau\approx
0.79$ is obtained for $\rho=\frac12\sqrt{9-\sqrt{73}}$. By the
$\lambda$-Lemma this remains true in all of $\W_1$. ~$\square$

\begin{thm} \label{T6.2}Any component $W$ of $\W_n$ $(n>1)$ is
simply connected and contains exactly one zero of the function
$Q_{n-2}^2-2$. The multiplier map is a degree-eight proper map
$W\longrightarrow\D$; it is ramified exactly at the zero $t_n$ of
$Q_{n-2}^2-2$, thus $\W_n$ has at most $2\deg Q_{n-2}
(t)=2(4^{n-1}-1)/3$ components.\end{thm}

{\bf Proof.} Theorems like that hold in any family with single
free critical orbit
$$c_0(t)\xrightarrow{(m_1)}c_1(t)\xrightarrow{(m_2)}\cdots
 \xrightarrow{(m_{k})}c_{k}(t)\xrightarrow{(1)}\cdots\quad(m_1>1, m_k>1)$$
and $\deg_{c_0(t)} f^\ell=m=m_1m_2\cdots m_k$ for $\ell\ge k$.
This $m$ is also the degree of the multiplier map
$\lambda:W\longrightarrow\D$. Details are worked out in
\cite{Ste3} and will be omitted. In our case we have $({\it
one~of})\pm\sqrt{2}\xrightarrow{(2)}0\xrightarrow{(4)}
 t\xrightarrow{(1)}\cdots$
and $m=8$. The number of hyperbolic components of order $n$ is at
most the number of {\it centres}, i.e.\  parameters such that
$f^n_t$ has a super-attracting periodic point of exact period $n$,
hence this number equals the number of solutions to the equation
$Q_n(t)=t$ (equivalently $Q_{n-2}(t)^2=2$) that are not solutions
to any equation $Q_{\ell}(t)=t$ with $\ell<n$. ~$\square$

\begin{thm} \label{T6.4}Suppose $t\in \W_n$ for some $n>1$. Then the
Julia set $\J_t$ is connected and locally connected. The Fatou set
consists of the basin $\A_t$ and its pre-images, and of the
domains $U_\nu$ $(0\le\nu <n)$ of the $($super-$)$attracting cycle
and their pre-images. Exactly one of the critical points $\pm\sqrt
2$ is contained in $\bigcup_{\nu=0}^{n-1}U_\nu$.
\end{thm}

\bigskip\includegraphics[viewport=0 0 368 153]{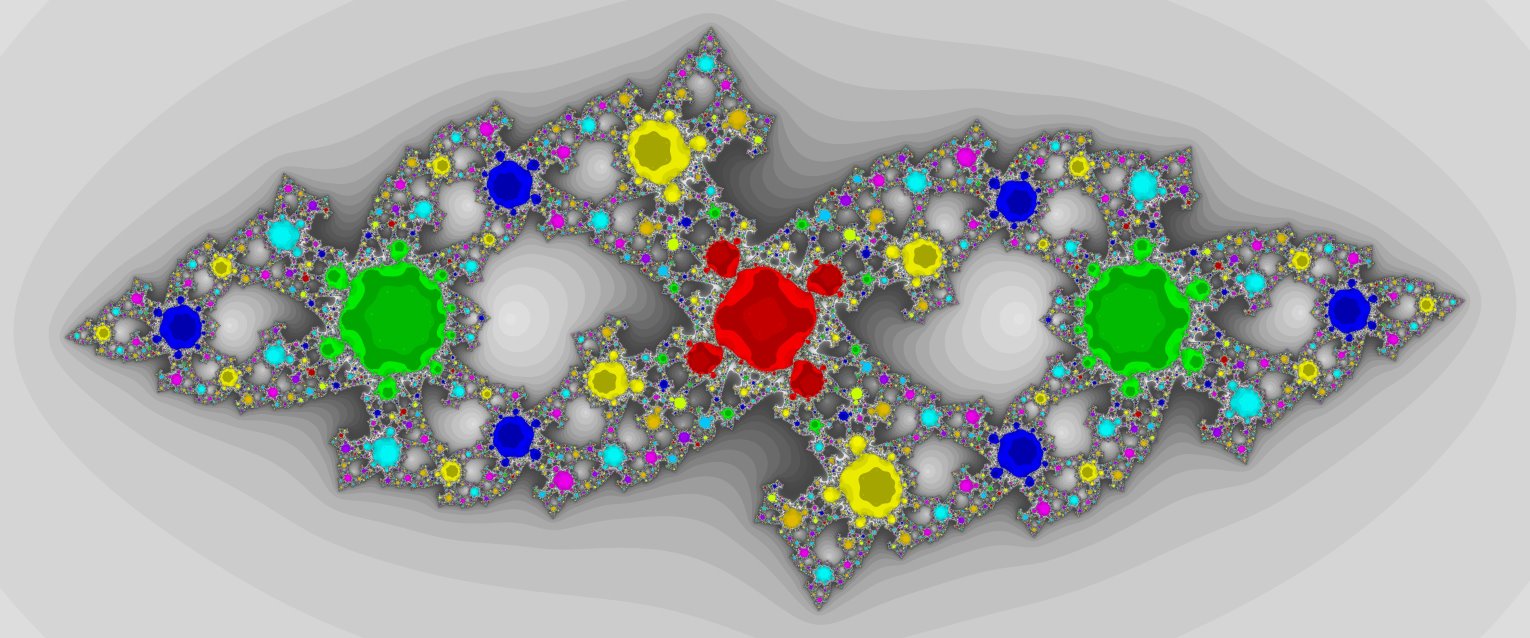}

\bigskip\includegraphics[viewport=0 0 368 184]{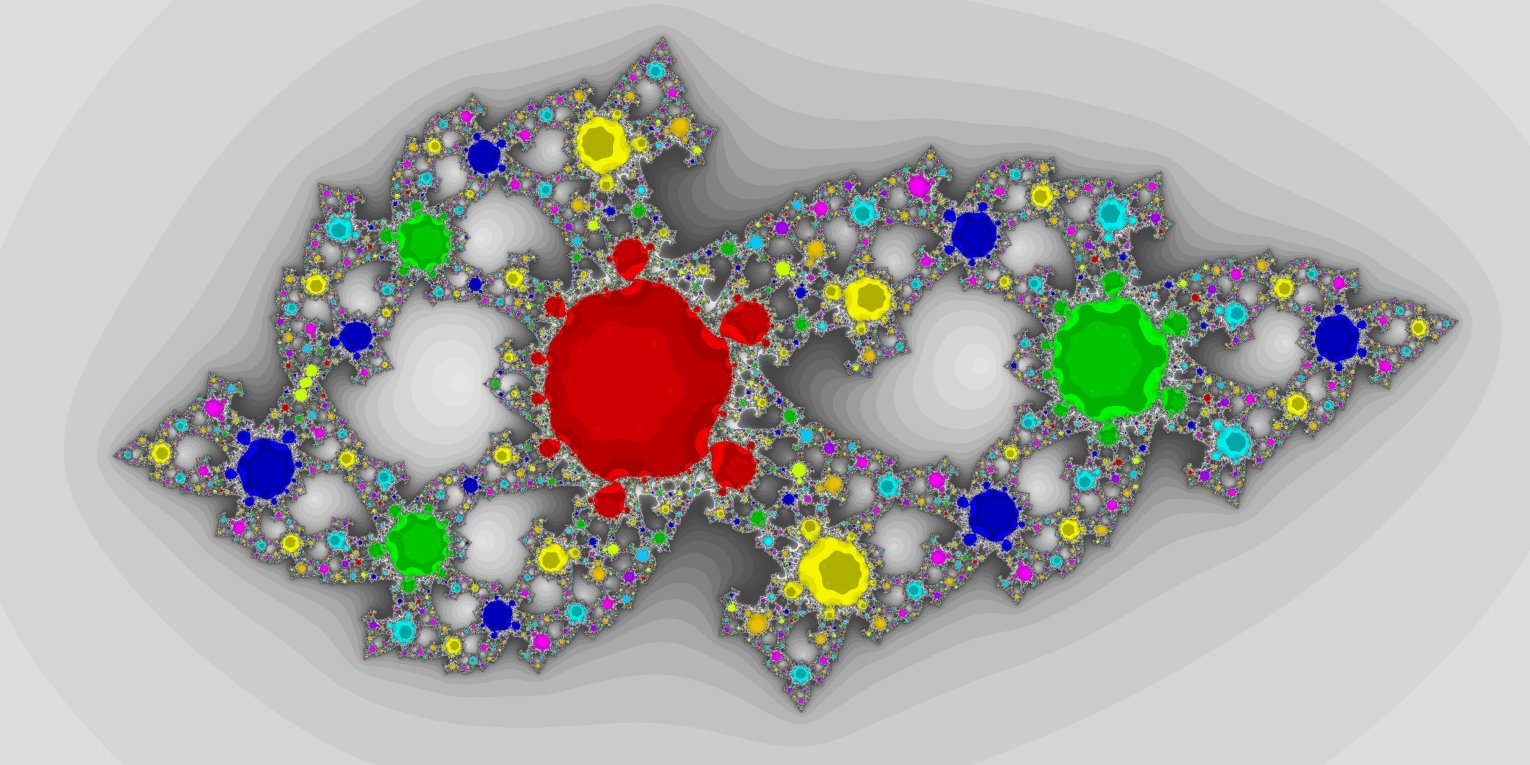}
\medskip

\begin{center}\small {\bf Figure 5.} $f_t$ ($t={1.33-0.133 i}$) as well as
its quasi-conjugate $\hat f_t(z)=f_t(\sqrt{z})^2$ has an
attracting $10$-cycle; the Julia set is connected and locally
connected.\end{center}

{\bf Proof.} We start with the {\it semi-conjugate} $\ds\hat
f_t(z)=(f_t(\sqrt{z}))^2=\frac{t^2}{16}\frac{(z-2)^4}{(z-1)^2}$
with free critical orbit
$2\xrightarrow{(4)}0\xrightarrow{(2)}t^2\xrightarrow{(1)}\cdots$.
Its super-attracting basin $\hat\A_t$ about $\infty$ as well as
the pre-images of $\hat\A_t$ of any order are simply connected.
Any attracting $n$-cycle $\{\hat U_0 ,\hat U_1 ,\ldots,\hat
U_{n-1}\}$ contains the critical point $z=2$, we may assume $0\in
\hat U_1$ and $2\in \hat U_0$, so that $\hat f_t:\hat
U_0\xrightarrow{(4)}\hat U_1\xrightarrow{(2)} \hat
U_2\xrightarrow{(1)}\cdots\xrightarrow{(1)}\hat U_{n-1}$ and $
\hat f_t^n:\hat U_0\xrightarrow{(8)}\hat U_0$ with one critical
point $(z=2)$ of order seven. Thus the domains $\hat U_\nu$ are
simply connected, and so are the domains $U_\nu$ $(\hat U_\nu$ is
the image of $U_\nu$ under the map $z\mapsto z^2$). From the
Riemann-Hurwitz formula applied to $f_t:U_0\longrightarrow U_1$
follows that $f_t|_{U_0}$ cannot have degree four with two
critical points $\pm\sqrt{2}$, hence has degree two, and $U_0$
contains only one of the critical points $\pm\sqrt{2}$. Thus
$$f_t:U_0\xrightarrow{(2)}U_1\xrightarrow{(4)}
 U_2\xrightarrow{(1)}\cdots\xrightarrow{(1)}U_{n-1}\quad{\rm and}\quad f_t^n:U_0\xrightarrow{(8)}U_0.$$
Since all Fatou components are simply connected, the Julia set is
connected, and is locally connected since $f_t$ is
hyperbolic.~$\square$

\medskip\small{\it Hye Gyong Jang}\hfill {\bf
E-mail:} \texttt{pptayang@co.chesin.com}\\
{\bf Address:} {\it Faculty of Mathematics and Mechanics,
University of Science Pyongyang, D.P.R. of Korea.}\medskip

{\it Norbert Steinmetz}\hfill {\bf
E-mail:} \texttt{stein@math.tu-dortmund.de}\\
{\bf Address:} {\it Fakult\"at f\"ur Mathematik, Technische
Universit\"at Dortmund, Germany.}\\

\end{document}